\documentclass{amsart}

\usepackage{color}

\newcommand\remove[1]{}

\newtheorem{theorem}{Theorem}[section]
\newtheorem{lemma}[theorem]{Lemma}

\newtheorem{cor}[theorem]{Corollary}

\theoremstyle{definition}
\newtheorem{definition}[theorem]{Definition}
\newtheorem{example}[theorem]{Example}

\theoremstyle{remark}

\numberwithin{equation}{section}

\newcommand{\Om}{\Omega}
\newcommand{\La}{\Lambda}

\newcommand{\s}{\sigma}

\def \N {\mathbb{N}}

\newcommand{\sw}{\mathsf{w}}

\newcommand{\1}{\mathbf{1}}

\DeclareMathOperator{\supp}{{supp}}

\newcommand{\lb}{\label}

\newcommand{\wtw}{if and only if}

\begin{document}

\title{Disjointness preserving operators on K\"{o}the-Banach spaces}

\author{Iryna Krasikova}

\address[Krasikova]{Zaporizhzhya National University\\
str.~Zhukovskoho 66, Zaporizhzhya, 69002 (Ukraine)}

\email{irynazpukr@gmail.com}

\author{Mikhail Popov}

\address[Popov]{Institute of Exact and Technical Sciences,
Pomeranian University in S{\l}upsk (Poland)}

\email{misham.popov@gmail.com}

\author{Beata Randrianantoanina}\thanks{The third named author was supported by the National Science Foundation under Grant No. 2620717}

\address[Randrianantoanina]{Miami University\\
Oxford, OH, USA}

\email{randrib@miamioh.edu}

\subjclass[2020]{Primary 46E30; Secondary 47B65; 46B15}

\keywords{K\"{o}the-Banach space; r.i. space;  rearrangement-invariant space; disjointness preserving operator; regular operator; narrow operator, $L_p(L_q)$-spaces}

\begin{abstract}
A well-known result of Abramovich says that if $E$, $F$ are rearrangement invariant Banach spaces   on a finite atomless measure space such that $E$ has an absolutely continuous norm and  $E\not\subseteq F$, then there are no nonzero disjointness preserving  bounded linear operators $T \colon E\to F$.

We obtain an analogue of this result in the setting of general  K\"{o}the-Banach spaces. Specifically, we identify a sufficient condition on  K\"{o}the-Banach spaces $E$, $F$ on finite atomless measure spaces   that implies that there are no nonzero disjointness preserving  bounded linear operators $T \colon E\to F$. Consequently, by a representation theorem for regular operators due to O.~Maslyuchenko, Mykhaylyuk and the second named author, this condition implies that every regular linear operator  $T \colon E\to F$ is narrow. As  an application, we obtain that if  $1\le \max\{p_1,p_2\}< \min\{ r_1, r_2\}\le \infty$, then
every regular linear operator $T \colon L_{p_1}(L_{p_2}) \to L_{r_1}(L_{r_2})$
 is narrow. The latter may be seen as an analogue  of the classical Pitt compactness theorem.
\end{abstract}

\maketitle

\section{Introduction}

Abramovich \cite{Ab} proved the following well-known theorem.

\begin{theorem}\lb{Abth}
Let $E$, $F$ be two rearrangement invariant Banach spaces on a finite atomless measure space and $E$ has an absolutely continuous norm. If there exists a nonzero disjointness preserving  bounded linear operator $T \colon E\to F$,
then $E\subseteq F$.
\end{theorem}

Here  $E\subseteq F$ means that the identity mapping from $E$ to $F$ is a bounded operator. Hence, in particular, since $E$, $F$ are rearrangement-invariant,
 there exists a constant $M<\infty$ such that for all sets $A,B\in \Sigma$, with $\mu(A)=\mu(B)$ we have $\|\1_B\|_F\le M \|\1_A\|_E$, that is,
\begin{equation}
 \frac{\|\mathbf{1}_A\|_E}{\|\mathbf{1}_B\|_F} \ge \frac{1}{M}.
\end{equation}

Thus, Abramovich's Theorem~\ref{Abth} immediately implies the following result.
\begin{theorem}\lb{Abcp}
Let $E$, $F$ be two rearrangement invariant Banach spaces on a finite  atomless measure space $(\Omega, \Sigma, \mu)$ and $E$ has an absolutely continuous norm.  If
\begin{equation} \label{eq:abcp}
\inf_{0<t\le 1} \sup \left\{ \frac{\|\mathbf{1}_A\|_E}{\|\mathbf{1}_B\|_F} : \, \frac{\mu(A)}{ \mu(\Omega)}
=\frac{\mu(B)}{ \mu(\Omega)}= t \right\} = 0,
\end{equation}
then there are no
nonzero  disjointness preserving  bounded linear operators $T \colon E\to F$.
\end{theorem}

 Theorem~\ref{Abcp} is slightly weaker than the contrapositive of  Theorem~\ref{Abth}, since 
  there exist
 non-isomorphic  rearrangement invariant  spaces whose norms coincide on all characteristic functions $1_A$, for $A\in\Sigma$, see
 \cite[Section~2.5]{BS} for full description of such spaces.

In the present note, we generalize Theorem~\ref{Abcp}  to  a  setting of not necessarily rearrangement invariant  K\"{o}the-Banach spaces.  As shown by Abramovich in \cite{Ab}, the assumptions that the measure space is atomless and that $E$ has an absolutely continuous norm are both necessary in Theorem~\ref{Abth}.
However we are able to completely remove the assumption that the space  $E$ is rearrangement invariant, and weaken the assumption on the space  $F$ and  require only that the averaging projection operators on  $F$ are uniformly bounded  (see \S\ref{sec:aver} for details). Specifically, we prove the following generalization of Theorem~\ref{Abcp}.

\begin{theorem}[see Theorems~\ref{cor:oajvn} and \ref{cor:ksjd2}] \label{cor:oajvn-intro}
Let $E,F$ be K\"{o}the-Banach spaces on finite atomless measure 
spaces $(\Omega_E, \Sigma_E, \mu)$ and $(\Omega_F, \Sigma_F, \nu)$, respectively, such that $E$ has an absolutely continuous norm and $F$ is averaging friendly. If
\begin{equation} \label{eq:sjnd}
\inf_{0<t\le 1} \sup \left\{ \frac{\|\mathbf{1}_A\|_E}{\|\mathbf{1}_B\|_F} : \,A\in  \Sigma_E, B\in \Sigma_F, {\text{ with }}  \frac{\mu(A)}{ \mu(\Omega_E)}=\frac{\nu(B)}{\nu(\Omega_F)} = t \right\} = 0,
\end{equation}
then there are no nonzero disjointness preserving  bounded linear operators $T \colon E\to F$, and every regular
linear operator $T \colon E\to F$ is narrow.
\end{theorem}

The second part of the conclusion follows from the first part combined with the representation theorem for regular operators due to O.~Maslyuchenko, Mykhaylyuk and the second named author \cite{MMP} (see \S\ref{sec:reg} below).

As one application of Theorem~\ref{cor:oajvn-intro} we obtain  that $\max\{p_1,p_2\}< \min\{ r_1, r_2\}$  then
  every regular linear  operator  $T \colon L_{p_1}(L_{p_2}) \to L_{r_1}(L_{r_2})$ is narrow, which may be seen as an analogue of the Pitt compactness theorem \cite{Pitt}, see Corollary~\ref{cor:LpLq}.  A similar result in the case of
   $L_p$-spaces was proved in \cite{MP4}.

Our methods are quite different from the methods of  \cite{Ab}. Our approach uses the theory of narrow operators.
The key ingredient of the proof of  Theorem~\ref{cor:oajvn-intro} is the fact that all continuous linear functionals on K\"{o}the-Banach spaces with an absolutely continuous norm are
 strictly narrow (see \S\ref{sec:narrow}).

\section{Preliminaries}

We follow standard notation of Banach space theory and Banach lattice theory, as may be found in \cite{ABu,BS,LTII,MN,PR,Schaefer}, where we also refer the reader for more background information and undefined terms.

\subsection{K\"{o}the-Banach spaces and  rearrangement invariant spaces}

Let $(\Omega, \Sigma, \mu)$ be a finite measure space. A Banach space $E$ of (equivalence classes of) measurable functions on $\Omega$ is called a \textit{K\"{o}the-Banach space} if the following conditions hold:

\begin{enumerate}
  \item[$(K_i)$] if $y \in E$ and $|x| \leq |y|$ then $x \in E$ and $\|x\| \leq \|y\|$;
  \item[$(K_{ii})$] $\mathbf{1}_\Omega \in E$;
  \item[$(K_{iii})$] $E \subseteq L_1(\mu)$.
\end{enumerate}

Let $E$ be a K\"{o}the-Banach space. By $(K_{iii})$, the inclusion embedding of $E$ to $L_1(\mu)$ is continuous and every element $y \in L_1(\mu)$ such that $x \cdot y \in L_1(\mu)$ for every $x \in E$, defines an element $g \in E^*$ by $g(x) = \int_\Omega x y \, d \mu$.

A K\"{o}the-Banach space $E$ is said to have an \textit{absolutely continuous norm}, if $\lim_{\mu(A) \to 0}  \|x \cdot \mathbf{1}_{A}  \| = 0$ for each $x \in E$. Some authors use the terminology that $E$ has an  \textit{order continuous norm}. By    \cite[Theorem~2.6]{ABu2} and \cite[Proposition~3.5]{BS}, these two notions are equivalent in  K\"{o}the-Banach spaces.

The K\"{o}the-Banach space  $L_p(\mu)$ has an absolutely continuous norm \wtw\ $1 \leq p < \infty$. 

For   $f \in L_1(\mu)$,  the \textit{distribution function} $d_{f}: [0,\infty) \to [0,\mu(\Omega)]$ of $f$ is defined by $d_{f}(a) = \mu \{\omega \in \Omega: \, f(\omega) > a \}$.  Functions $f$ and $g$ in $L_1(\mu)$ are called \textit{equimeasurable} if they have the same distribution function, that is, $d_{f}(a)=d_{g}(a)$ for all $a \in \mathbb R$.  Functions $f$ and $g$ are called \textit{equimeasurable in modulus} if their moduli have the same distribution function, that is, $d_{|f|}(a)=d_{|g|}(a)$ for all $a\ge 0$.

A  K\"{o}the-Banach space $E$ on a $\s$-finite measure space is called {\it rearrangement invariant}  (r.i., for short) if
\begin{enumerate}
\item[$(K_{iv})$] if $x\in E$ and $d_{|y|}=d_{|x|}$ then $y\in E$ and $\|y\|_E=\|x\|_E$.
\end{enumerate}

Rearrangement invariant spaces have been extensively studied in the literature. For example,
all  spaces $L_p(\mu)$ with $1 \leq p \leq \infty$, Orlicz spaces, and Lorentz spaces are r.i. spaces. For more information on r.i. spaces see \cite{BS,LTII}.
Examples of  K\"{o}the-Banach spaces that are not rearrangement invariant include weigthed $L_p$-spaces, weighted Orlicz spaces, spaces of the form $L_p(L_q)$ (considered as a space of functions on the product measure space), or $E(F)$, where $E$, $F$ are
K\"{o}the-Banach spaces, excluding only the case when $E=F=L_p$
for some $p$ (since $L_p(L_p)=L_p $ is an r.i. space), and many others.

\subsection{Averaging friendly K\"{o}the-Banach spaces}\lb{sec:aver}

\begin{definition} \label{def:ue}
 We say that a K\"{o}the-Banach space $F$ on an 
   measure space $(\Omega, \Sigma, \nu)$ is \emph{$K$-averaging friendly}, where $K \in [1, +\infty)$, if for every $B \in \Sigma$ with $\nu(B) > 0$, the averaging projection operator $P_{B} \colon F \to F$ defined by
\begin{equation} \label{eq:eoncd}
P_{B} (g) := \left( \frac{1}{\nu(B)} \int_B g \, d \nu \right) \cdot \mathbf{1}_B,\ \ \ \ \ \ (g \in F)
\end{equation}
satisfies $\|P_{B}\| \le K$. We say that  $F$ is  \emph{averaging friendly}, if $F$ is $K$-averaging  friendly for some constant $K<\infty$.
\end{definition}

It is well-known (and easy to see) that every r.i. space  on an atomless measure space $(\Omega, \Sigma, \nu)$ is averaging friendly with constant $K = 1$ (see \cite[Theorem~2.a.4]{LTII}).
However, not every K\"{o}the-Banach space  is averaging friendly, even if it has an absolutely continuous norm.

\begin{example}[A K\"{o}the-Banach space that is  not averaging friendly] \label{ex:jjjU}

Let $\sw$ be a weight on $[0,1]$ defined by $\sw:= \sum_{n=0}^\infty 2^{-n} \mathbf{1}_{[2^{-(n+1)},2^{-n}]}$, and let
$L_{1,\sw}[0,1]$ be the weighted $L_1$-space, that is $L_{1,\sw}[0,1]=(L_1[0,1],\|\cdot\|_{1,\sw})$,  where
\begin{equation} \label{eq:kdsjj757}
\|x\|_{1,\sw} :=  \int_{[0,1]} \sw |x| \, d \mu=\sum_{n=0}^\infty 2^{-n} \int_{\left[\frac{1}{2^{n+1}}, \frac{1}{2^{n}} \right]} |x| \, d \mu.
\end{equation}
Then $L_{1,\sw}[0,1]$ is not averaging friendly.
\end{example}

\begin{proof}
Indeed, for $x_m := \mathbf{1}_{\left[\frac{1}{2^{m+1}}, \frac{1}{2^{m}} \right]}$ and $B = [0,1]$ one has
$$
\frac{\|P_B (x_m)\|_{1,\sw} }{\|x_m\|_{1,\sw} } = \frac{\frac23 \cdot 2^{-m-1}}{2^{-2m-1}} = \frac{2^{m+1}}{3}
$$
for all $m \in \mathbb N$, proving that the linear operator $P_{B}$ is unbounded on $F$.
\end{proof}

Using the idea of Example~\ref{ex:jjjU} and considering weights $\sw_M$ defined by  finite sums in place of the infinite sum, for any given $M > 1$ one can construct an averaging friendly K\"{o}the-Banach space on $[0,1]$, which is not averaging friendly with any constant $K < M$.

\subsection{The fundamental function}

Let $E$ be an r.i. space on $[0,1]$. The function
$$
\phi_E(t) := \|\mathbf{1}_{[0,t]}\|_E,\, \ \ \   (t \in (0,1]) \,\,\, 
$$
 is called the \emph{fundamental function of $E$}.

We are also interested in the following function:
$$
\psi_E(t) := \frac{\phi_E(t)}{t}, \, \ \ \ \ \   (t \in (0,1]). \,\,\,
$$

It is known that for all r.i. spaces $E$, the fundamental function $\phi_E$ is increasing, and $\psi_E$ is decreasing, see \cite[p.~67, Corollary~5.3]{BS}.

We introduce an analogue of the fundamental function for K\"{o}the-Banach spaces.

\begin{definition}
Let $E$ be a K\"{o}the-Banach space on a finite 
measure space $(\Omega, \Sigma, \mu)$.
The function $\phi_E \colon \Sigma \to [0, +\infty)$ defined by
$$
\phi_E (A) = \|\mathbf{1}_A\|_E, \,\, \ \ (A \in \Sigma),
$$
is called the \emph{fundamental function} of $E$. We also define the function $\psi_E \colon \Sigma^+ \to (0, +\infty)$ by setting
$$
\psi_E(A) = \frac{\phi_E (A)}{\mu(A)}, \, \,\ \  (A \in \Sigma^+ := \{A \in \Sigma: \, \mu(A) > 0\}).
$$
\end{definition}

Since $E$ is a  K\"{o}the-Banach space, the fundamental function is increasing, that is, for all  $A,B \in \Sigma$ with $A \subseteq B$, we have $\phi_E (A) \le \phi_E (B)$.

\begin{lemma} \label{le:Bjdnc}
Let $K<\infty$ and $F$ be a $K$-averaging friendly K\"{o}the-Banach space on a finite atomless measure space $(\Omega_F, \Sigma_F, \nu)$. Then for all $B',B \in \Sigma^+$,
\begin{equation*} B' \subseteq B \implies \psi_F (B) \le K \psi_F (B').
\end{equation*}
\end{lemma}

\begin{proof}
Since $ B' \subseteq B$ we have
$$P_B(\mathbf{1}_{B'})= \frac{\nu(B')}{\nu(B)}\mathbf{1}_{B}.$$ Thus $\frac{\nu(B')}{\nu(B)}\|\mathbf{1}_{B}\|\le K\|\mathbf{1}_{B'}\|$, and hence $\psi_F (B) \le K \psi_F (B')$.
\end{proof}

\subsection{Disjointness preserving operators}

Let $E$ be a  K\"{o}the-Banach space on $(\Om,\Sigma,\mu)$. We say that two elements $x, y$ of $E$ are {\it disjoint}, denoted  $x \perp y$, if $|x|\wedge |y|=0$.

We write $e = x_1 \sqcup \ldots \sqcup x_n$ to express that $e = x_1 + \ldots + x_n$ and $x_i \perp x_j$ for all $i \neq j$.
 Analogously, we use the notation $A = A_1 \sqcup \ldots \sqcup A_n$ for sets $A, A_i \in \Sigma$ to express that $A = A_1 \cup \ldots \cup A_n$ and $A_i \cap A_j = \emptyset$ for $i \neq j$.

Let $E$ and $F$ be K\"{o}the-Banach spaces. A bounded linear operator $T \colon E \to F$ is called \emph{disjointness preserving}, if for every $x,y \in E$ the condition $x \perp y$ implies $Tx \perp Ty$. It is easily seen that if $T \colon E \to F$ is disjointness preserving then $e = x_1 \sqcup \ldots \sqcup x_n$ for elements of $E$ implies that $Te = Tx_1 \sqcup \ldots \sqcup Tx_n$.

\subsection{Narrow and strictly narrow operators} \lb{sec:narrow}

Let $E$ be a K\"{o}the-Banach space on a finite atomless measure space $(\Omega, \Sigma, \nu)$ and $X$ be a Banach space. Following \cite[Definition~1.5]{PR}, cf. also  \cite[p.~53]{PP},  we say that a linear operator 
$T \colon E \to X$ is:
\begin{itemize}
\item \emph{narrow}, if for every set $A \in \Sigma$ and every $\varepsilon > 0$ there exists a decomposition $A = B \sqcup C$ such that $\nu(B) = \nu(C)$ and $\|T\mathbf{1}_B - T\mathbf{1}_C\| < \varepsilon$;
  \item \emph{strictly narrow}, if for every  $A \in \Sigma$  there exists a decomposition $A = B \sqcup C$ such that $\nu(B) = \nu(C)$ and $T\mathbf{1}_B = T\mathbf{1}_C$.
\end{itemize}

By the Liapunov convexity theorem \cite[Theorem~2.c.9]{LTII}, the range of every atomless $\sigma$-additive vector measure on a measurable space $(\Omega, \Sigma)$ taking values in a finite-dimensional normed space is convex. As a consequence, we obtain the following result about finite rank operators.

\begin{theorem} \label{th:jadc}
Let $E$ be a K\"{o}the-Banach space on a finite atomless measure space $(\Omega, \Sigma, \mu)$ with an absolutely continuous norm, and $X$ be  a finite-dimensional normed space. Then every bounded linear  operator $T \colon E\to X$  is strictly narrow.
\end{theorem}

\begin{proof}
Indeed, every bounded linear  operator $T \colon E\to X$  generates a vector measure $\nu_T \colon \Sigma \to X \oplus \mathbb K$, where $\mathbb K \in \{\mathbb R, \mathbb C\}$, defined by  $\nu_T (A) = \bigl(T \mathbf{1}_A, \mu(A) \bigr)$, $(A \in \Sigma)$. The measure$\nu_T$ is $\sigma$-additive by the absolute continuity of the norm on $E$. By the Liapunov convexity theorem, since $X \oplus \mathbb K$ is finite dimensional, the range of $\nu_T$ is convex. Hence, for every $A \in \Sigma$, the element $1/2\big(\nu_T(\emptyset)+\nu_T(A)\big)=\bigl(T \mathbf{1}_A/2, \mu(A)/2 \bigr)$
belongs to  the range  of $\nu_T$, which exactly means the strict narrowness of $T$.
\end{proof}

Theorem~\ref{th:jadc} is not new - it can also be obtained as a direct corollary of  \cite[Theorem~4.11]{MPPS}.
The assumption of the   absolute continuity of the norm  in Theorem~\ref{th:jadc} is essential: there exist nonnarrow continuous linear functionals on $L_\infty$, see e.g. \cite[Example~10.12]{PR}.

\section{Main results}

In this section we prove  the first part of  Theorem~\ref{cor:oajvn-intro}, which we restate here.

\begin{theorem} \label{cor:oajvn}
Let $E,F$ be K\"{o}the-Banach spaces on  finite  atomless measure 
spaces $(\Omega_E, \Sigma_E, \mu)$ and $(\Omega_F, \Sigma_F, \nu)$, respectively, such that $E$ has an absolutely continuous norm and $F$ is averaging friendly. If
\begin{equation} \label{eq:sjnd}
\inf_{0<t\le 1} \sup \left\{ \frac{\|\mathbf{1}_A\|_E}{\|\mathbf{1}_B\|_F} : \, A\in  \Sigma_E, B\in \Sigma_F, {\text{ with }} \frac{\mu(A)}{ \mu(\Omega_E)}=\frac{\nu(B)}{\nu(\Omega_F)} = t \right\} = 0,
\end{equation}
then there are no nonzero disjointness preserving bounded linear operators
from $E$ to $F$.
\end{theorem}

\begin{proof}
Assume, for contradiction, that there exists a nonzero disjointness preserving bounded linear  operator
$T:E\to F$.
We will show that then  there exists a constant
$M < \infty$ such that for every $t \in (0,1]$, there exist sets $A \in \Sigma_E$ and $B \in \Sigma_F$ such that
$ \mu(A)/\mu(\Omega_E)=\nu(B)/\nu(\Omega_F)=t$ and
\begin{equation}\lb{bigratio}
 \frac{\|\mathbf{1}_A\|_E}{\|\mathbf{1}_B\|_F} \ge \frac{1}{M},
\end{equation}
which contradicts \eqref{eq:sjnd} and thus will end the proof of the theorem.

Let  $T \colon E\to F$ be a nonzero disjointness preserving bounded linear  operator.
First we show that $y := T \mathbf{1}_{\Omega_E} \neq 0$. Indeed, since $T$ preserves disjointness, for every $C \in \Sigma_E$ one has $T \mathbf{1}_C \sqsubseteq y$, and the assumption $y = 0$ implies $T \mathbf{1}_C = 0$ for all $C \in \Sigma_E$. Thus $T|_S = 0$, where $S$ is the set of all simple (= finite valued) functions. Since the norm of $E$ is absolutely continuous, $S$ is dense in $E$, see \cite[Proposition~2.10]{PR} (there exists a K\"{o}the-Banach space on $[0,1]$, in which the set of all simple functions is not dense, see \cite{PS}). Thus, $T = 0$, a contradiction. So $y \neq 0$.

Let  $s := {\rm sign} \, y$. We define  a functional $f \in E^*$ by setting
$$
f(x) = \int_{\Omega_F} (Tx) \cdot s \, d \nu, \ \ \ \ \  (x \in E).
$$

Observe that
$$
f(\mathbf{1}_{\Omega_E}) = \int_{\Omega_F} y \cdot s \, d \nu = \int_{\Omega_F} |y| \, d \nu = \|y\|_{L_1}> 0.
$$

Fix any $t\in(0,1]$, and let  $n \in \mathbb N$, $m\in\{1,\dots,2^n\}$, be such that  $t/2\le \frac{m}{2^n}\le t$.

By Theorem~\ref{th:jadc}, $f$ is strictly narrow.  Thus, by  definition, we can choose a decomposition $\Omega_E = A_1 \sqcup \ldots \sqcup A_{2^n}$ so that $\mu (A_k) = 2^{-n} \mu (\Omega_E)$ and $f(\mathbf{1}_{A_k}) = 2^{-n} f(\mathbf{1}_{\Omega_E})$, for all $k \in \{1, \ldots, 2^n\}$.

For $k \in \{1, \ldots, 2^n\}$, set $x_k := \mathbf{1}_{A_k}$, $y_k := T x_k$, and $B_k := \supp y_k$.
Note that $f(x_k)=2^{-n} f(\mathbf{1}_{\Omega_E})=2^{-n}\|y\|_{L_1}> 0$. Thus $Tx_k\ne 0$, which implies that
 $\nu(B_k)>0$ for all $k$. Since the operator $T$ is disjointness preserving, the sets $\{B_k\}_{k=1}^{2^n}$ are mutually disjoint. Therefore,
 \[\nu(\supp y)=\sum_{k=1}^{2^n} \nu(B_k).\]

Let $\s$ be a permutation of $\{1,\dots,2^n\}$, such that the sequence  $\{\nu(B_{\s(k)})\}_{k=1}^{2^n}$ is
nondecreasing. Given any  $ m\in\{1,\dots,2^n\}$, define $\La:=\{\s(1),\dots,\s(m)\}\subseteq \{1,\dots,2^n\}$,
$A_\La:=\bigsqcup_{k\in\La} A_k$, $x_\La := \mathbf{1}_{A}$, $y_\La:= T x_\La$, and $B_\La:= \supp y_\La$.

Since  $|\La|=m$, $\mu(A_\La)= \frac{m}{2^n} \mu (\Omega_E)\in[\frac{t}{2}\mu (\Omega_E), t \mu (\Omega_E)]$, and, since the sequence  $\{\nu(B_{\s(k)})\}_{k=1}^{2^n}$ is
nondecreasing, we have
\[0<\nu(B_\La)=\sum_{k=1}^{m}  \nu(B_{\s(k)})\le  \frac{m}{2^n} \sum_{k=1}^{2^n}  \nu(B_{\s(k)})\le
 t \nu(\supp y)\le  t \nu(\Om_F).\]

Let $K<\infty$ be such that $F$ is $K$-averaging friendly. Then
we obtain
\begin{equation} \label{eq:jsah}
\begin{split}
\frac{\|\mathbf{1}_{B_\La}\|_F}{\nu(B_\La)}
 &= \frac{1}{f(x_\La)} \left\| \left(\frac{1}{\nu(B_\La)} \int_{\Omega_E} (Tx_\La) \cdot s \, d \nu \right) \cdot
  \mathbf{1}_{B_\La} \right\|_F \\
&= \frac{2^n}{m f(\mathbf{1}_{\Omega_E})} \left|\frac{1}{\nu(B_\La)}\int_{B_k} (Tx_\La) \cdot s \, d \nu \right|
 \cdot \|\mathbf{1}_{B_\La}\|_F \\
 &\le \frac{2}{t \|y\|_{L_1}} \left(\frac{1}{\nu(B_\La)}\int_{B_\La} |Tx_\La|  \, d \nu \right)
 \cdot \|\mathbf{1}_{B_\La}\|_F \\
  &= \frac{2}{t \|y\|_{L_1}} \left\| P_{B_\La}\left( |Tx_\La|\right) \right\|_F\\
&\le \frac{2K}{t \|y\|_{L_1}} \big\| \, |Tx_\La|\, \big\|_F
= \frac{2K}{t \|y\|_{L_1}} \big\|Tx_\La \big\|_F
\le \frac{2K \|T\|}{t \|y\|_{L_1}}\, \|\mathbf{1}_{A_\La}\|_E.
\end{split}
\end{equation}

Let  $A \in \Sigma_E$   be such that $A_\La \subseteq A$ and $\mu(A) = t \mu(\Omega_E)$. Then
 $\|\mathbf{1}_{A_\La}\|_E\le \|\mathbf{1}_A\|_E$.
Let and $B \in \Sigma_F$ be  such that $B_\La \subseteq B$ and $\nu(B) = t \nu(\Omega_F)$. Then, by
Lemma~\ref{le:Bjdnc},
\begin{equation}\label{eq:t}
\begin{split}
\frac{\|\mathbf{1}_{B}\|_F}{t \nu(\Omega_F)} &= \frac{\|\mathbf{1}_{B}\|_F}{\nu(B)} = \psi_F(B) \le K \psi_F(B_\La)\underset{\eqref{eq:jsah}}{\le} \frac{2K^2 \|T\|}{t \|y\|_{L_1}} \|\mathbf{1}_{A_\La}\|_E\\
&\le  \frac{2K^2 \|T\|}{t \|y\|_{L_1}} \|\mathbf{1}_A\|_E.
\end{split}
\end{equation}

 Set
 $$M =: \frac{2K^2\|T\| \, \nu(\Omega_F)}{\|y\|_{L_1}}.$$
Then \eqref{eq:t} implies that $\|\|\mathbf{1}_A\|_E/\mathbf{1}_B\|_F \ge 1/M $, which ends the proof of \eqref{bigratio}, and the resulting contradition with \eqref{eq:sjnd} ends the proof of
the theorem.
\end{proof}

The following example shows that the atomlessness assumption on measure spaces in  Theorem~\ref{cor:oajvn} is essential.

\begin{example} \label{ex:lekfv}
Set $\Omega = \mathbb N$, $\Sigma = \mathcal{P}(\mathbb N)$ (where  $\mathcal{P}(\mathbb N)$ denotes the collection of all subsets of $\N$) and define the probability measure $\mu$ on $\Sigma$ by setting
$$
\mu(A) = \sum_{n \in A} 2^{-n}, \,\, \ \   (A \in \Sigma).
$$
Let $E = L_1(\mu)$ and $F = L_\infty(\mu)$. Then $E$ and $F$ are $1$-averaging friendly K\"{o}the-Banach spaces, $E$ has an absolutely continuous norm, and the following nonzero operator $T:E\to F$ is bounded ($\|T\|\le 1$) and preserves disjointness:
$$
T (\xi_n)_{n \in \mathbb N} = ( 2^{-n} \xi_n)_{n \in \mathbb N}, \,\, (\xi_n)_{n \in \mathbb N} \in L_1(\mu).
$$

However, for all $A, B\in \Sigma$, we have $\mu(A)=\mu(B)$ if and only if $A=B$. Moreover, for every $A\in \Sigma$,
$\|\mathbf{1}_A\|_E = \mu(A)$ and $\|\mathbf{1}_A\|_F = 1$, which implies that \eqref{eq:sjnd} is satisfied.
\end{example}

\section{Applications to regular operators between K\"{o}the-Banach spaces}\lb{sec:reg}

Let $E$ and $F$  be Riesz spaces.
A linear operator $T \colon E \to F$ is called:
\begin{itemize}
  \item \emph{positive} (briefly, $T \ge 0$ or $T\in L^+(E,F)$), if $T(E^+) \subseteq F^+$;
  \item \emph{regular}  ($T\in L_r(E,F)$) if $T$ is a difference of two positive operators;
  \item \emph{order bounded}    ($T\in L_b(E,F)$) if  $T$ sends order bounded subsets of $E$ to order bounded subsets of $F$;
  \item a \emph{lattice homomorphism} if $T(x \vee y) = Tx \vee Ty$ for all $x,y \in E$.
\end{itemize}

Every disjointness preserving operator is regular, see e.g. \cite[Theorems~2.40 and 1.18]{ABu}.

We denote the space of all  bounded linear operators $T \colon E \to F$ by $ \mathcal L(E,F)$.
If $E, F$ are Banach lattices then $L_r(E,F) \subseteq  \mathcal L(E,F)$ \cite[Proposition~1.3.5]{MN}.

If $F$ is Dedekind complete then, by the Riesz-Kantorovitch Theorem (see e.g. \cite[Theorem~1.18]{ABu}), $L_r(E,F) = L_b(E,F)$, and
the band in $L_r(E,F)$ generated by all disjointness preserving operators from $E$ to $F$ (that is, by all lattice  homomorphisms) is denoted by $\mathcal{H}(E,F)$. Following \cite{PR},  operators from $\mathcal{H}(E,F)$ are called {\it pseudo embeddings}.

By the Riesz theorem \cite[Theorem~1.42]{ABu}, if $F$ is Dedekind complete then we have the decomposition of $L_r(E,F)$ into the direct sum of two bands
\begin{equation}\label{decomposition}
L_r(E,F) = \mathcal{H}(E,F) \oplus \mathcal{D}(E,F).
\end{equation}

Using terminology of Weis \cite{W} and Huijsmans-de Pagter \cite{HP}, elements of the disjoint complement $\mathcal{D}(E,F) := \mathcal{H}(E,F)^d$ are called  \emph{diffuse operators}.

For a long time, the decomposition \eqref{decomposition} was the core of many representation theorems by different authors (Huijsmans, Kalton, Kwapien, O.~Maslyuchenko, de~Pagter, Pliev, Rosenthal, Weis and the second-named author), the purpose of which is description of the summands. We refer the reader to \cite{PP24} for the history of representation theorems for regular operators.

We recall a consequence of the general description of  diffuse operators from \cite{PP24}.

\begin{theorem}[{cf. \cite[Theorem~B]{PP24}}]  \label{cor:ksjd}
Let $E$ and $F$ be  K\"{o}the-Banach spaces on atomless  measure
spaces. If there are no nonzero disjointness preserving  operators from $E$ to $F$, then every regular operator from $E$ to $F$ is narrow.
\end{theorem}

Combining Theorem~\ref{cor:ksjd} and  Corollary~\ref{cor:oajvn}, we obtain the following result.

\begin{theorem} \label{cor:ksjd2}
Let $E,F$ be K\"{o}the-Banach spaces on finite atomless  measure  
spaces $(\Omega_E, \Sigma_E, \mu)$ and $(\Omega_F, \Sigma_F, \nu)$ respectively. Assume that $E$ has an absolutely continuous norm and $F$ is averaging friendly. If \eqref{eq:sjnd} holds then every regular operator from $E$ to $F$  is narrow.
\end{theorem}

Note that Theorem~\ref{cor:ksjd2} for the partial case of $L_p$-spaces was mentioned in a recent paper by Mykhaylyuk and the second-named author \cite{MP4}, where they proved, among other results, that  if  $1 \le p < r \le \infty$, then
 every regular operator $T \colon L_p(\nu) \to L_r(\nu)$ is narrow.

Theorem~\ref{cor:ksjd2} can  be applied for example to spaces of the form $L_p(L_r)$ considered as
K\"{o}the-Banach spaces on the product $[0,1]\times[0,1]$ with the product Lebesgue measure. Thus we obtain:

\begin{cor} \label{cor:LpLq}
(a) If $1\le p< \min\{ r_1, r_2\}\le \infty$, then
every regular operator $T \colon L_{p} \to L_{r_1}(L_{r_2})$
 is narrow.

 (b) If $1\le \max\{p_1,p_2\}< r\le \infty$, then
every regular operator $T \colon L_{p_1}(L_{p_2}) \to L_r$
 is narrow.

(c) If $1\le \max\{p_1,p_2\}< \min\{ r_1, r_2\}\le \infty$, then
every regular operator $T \colon L_{p_1}(L_{p_2}) \to L_{r_1}(L_{r_2})$
 is narrow.
\end{cor}

\begin{proof}
Indeed, by Jensen's inquality, for every $1 \le q,r \le \infty$ and every measurable set
$A\subseteq [0,1]\times[0,1]$ we have
\[\mu(A)^{1/\min\{q,r\}}\le \| \mathbf{1}_A\|_{L_q(L_r)}\le \mu(A)^{1/\max\{q,r\}}.
\]
Thus, in each of the cases (a), (b), and (c),  condition  \eqref{eq:sjnd} is satisfied, and the conclusion follows by Theorem~\ref{cor:ksjd2}.
\end{proof}

It is easy to see that analogous statements are also true for iterated spaces of the form $L_{p_1}(L_{p_2}(...(L_{p_n})...))$.

\vspace{0,3 cm}

\centerline{\textbf{Data Availability Statement}}

Data sharing not applicable to this article as no datasets were generated or analyzed during the preparation of the paper.

\vspace{0,3 cm}

\section*{Declarations}

\subsection*{Conflict of interest}

On behalf of all authors, the corresponding author states that there is no conflict of
interest.

\bibliographystyle{amsplain}

\begin{thebibliography}{10}

\bibitem{Ab} Y.\,A.~Abramovich, \emph{Operators preserving disjointness on rearrangement invariant spaces}, Pacific J. Math. 148 (1991), no.~2, 201--209.

\bibitem{ABu2} C.\,D.~Aliprantis, O.~Burkinshaw. \emph{Locally solid Riesz spaces}. Acad. Press, New York, (1978).

\bibitem{ABu} C.\,D.~Aliprantis, O.~Burkinshaw, \emph{Positive Operators}, Springer, Dordrecht. (2006).

\bibitem{BS} C.~Bennett, R.~Sharpley, \emph{Interpolation of Operators}. Pure and Applied Mathematics, Vol. 129. Boston: Academic Press, 1988.

\bibitem{HP} C.\,B.~Huijsmans, B.~de Pagter, \emph{Disjointness preserving and diffuse operators}, Compositio Math., {\bf 79} (1991), no~3, 351--374.

\bibitem{LTII} J.~Lindenstrauss, L.~Tzafriri, \emph{Classical Banach spaces, Vol. 2, Function spaces}, Springer--Verlag, Berlin--Heidelberg--New York (1979).

\bibitem{MMP}  O.\,V.~Maslyuchenko, V.\,V.~Mykhaylyuk and M.\,M.~Popov, \emph{A lattice approach to narrow operators}, Positivity, 13 (3) (2009), 459-495. https://doi.org/10.1007/s11117-008-2193-z

\bibitem{MPPS} V.~Mykhaylyuk, M.~Pliev, M.~Popov, O.~Sobchuk, \emph{Dividing measures and narrow operators}, Stud. Math. 231 (2) (2015), 97--116. https://doi.org/10.4064/sm7878-2-2016

\bibitem{MP4} V. Mykhaylyuk, M. Popov. Some results on regular operators on $L_p$. Preprint.

\bibitem{MN} P.~Meyer-Nieberg, \textit{Banach Lattices}, Graduate Texts in Math., 233, Springer-Verlag, Berlin-Heidelberg-New York, 1991.

\bibitem{Pitt} H.\,R.~Pitt, \emph{A note on bilinear forms}, J. London Math. Soc., 11 (1936), pp.~174--180.

\bibitem{PP} A.~M.~Plichko, M.~M.~Popov, \emph{Symmetric function spaces on atomless probability spaces}, Dissertationes Math. (Rozprawy Mat.) 306 (1990), pp.~1--85.

\bibitem{PP24} M.~Pliev, M.~Popov. \emph{Representation theorems for regular linear and orthogonally additive operators}. Math. Nachr. 297, no 1 (2024), 126--143. https://doi.org/10.1002/mana.202200129

\bibitem{PR} M.~Popov, B.~Randrianantoanina, \emph{Narrow Operators on Function Spaces and Vector Lattices}, De Gruyter Studies in Mathematics 45, De Gruyter, 2013. https://doi.org/10.1515/9783110263343

\bibitem{PS} M.~Popov, O.~Sobchuk, \emph{On the ``function'' and ``lattice'' definitions of a narrow operator}, Positivity 22 (2018), 59--62. https://doi.org/10.1007/s11117-017-0497-6

\bibitem{Schaefer} Schaefer~H.~H. Banach Lattices and Positive Operators. Springer-Verlag. Berlin-Heidelberg-New York, 1974.

\bibitem{W} L.~Weis, \emph{On the representation of order continuous operators by random measures,} Trans. Amer. Math. Soc. \textbf{285} (1984), no.~2, 535--563.
\end{thebibliography}

\end{document}